\documentclass[11pt]{amsart}
\usepackage[T2A]{fontenc}
\usepackage[cp1251]{inputenc}
\usepackage[russian]{babel}
\usepackage{amsmath}
\usepackage{amssymb}
\usepackage{amscd}
\begin{document}
	
	\title{Invariants of singular germs of real-analytic sets}
	
	\author{M.\,A.~Stepanova}
	\address{Steklov Mathematical Institute of Russian Academy of Sciences, Moscow, Russia}
	\email{step\_masha@mail.ru}
	
	\thanks{This work was supported by the Russian Science Foundation under grant no 24-11-00196, https://rscf.ru/en/project/24-11-00196/}

	\renewcommand{\abstractname}{Abstract} 
	\begin{abstract} New biholomorphic invariants of a germ of a real-analytic set at a singular point (either an RC-singular point or a singular point of a smooth structure) are found: the Bloom-Graham singular type and the singular l-type. Their construction is based on the Bloom-Graham type and l-type constructions, respectively.
		
		A series of examples is considered.

	
	
	
	\end{abstract}
	\maketitle
	
	\section{Introduction}

	\vspace{3ex}
	The Bloom–Graham type (see \cite{1}) is an important biholomorphic invariant of a CR manifold germ. To construct it, a representative of a germ of $M$, defined in a neighborhood $U$ of the center $p$ of the germ, is associated with a chain of distributions $D_{\nu}$ of subspaces of complexifications of tangent spaces to the manifold (see Section 2). This gives a chain of nested subspaces of the complexification of the tangent space at the center of the germ (the values $D_{\nu}(p)$ of distributions at the center $p$ of the germ). The type is constructed from the sequence of dimensions of these subspaces.

	However, along with CR manifolds, it is natural to consider a wider class of sets whose elements admit singularities. Namely, let $M$ be locally defined as the set of common zeros of real-analytic functions in a complex space, and we require that, at a generic point (outside some proper analytic subset $\Sigma$), $M$ is a CR manifold of fixed CR type $(n,k)$. It is natural to call such sets \textit{CR-regular sets}, and to call the points at which $M$ is a CR manifold \textit{CR-regular} points.
	All our considerations and constructions are local and concern germs, and to construct auxiliary objects, we choose some representative $M$ of the germ $M_{p}$ at the point $p$.
	We will consider germs of locally connected (i.e., $M$ must be connected in some neighborhood of a singular point) and irreducible CR-regular sets (considering reducible sets reduces to considering their irreducible components). In the complement $\Sigma$ of a regular set ${\rm reg} \, M$, the set $M$ has singularities: either RC-singular points (i.e., points at which the dimension of the complex tangent space is greater than at the generic point) or singular points of a smooth structure. Our definition of a CR-regular set allows us to consider these two types of singularities from a unified perspective.
	
	The question arises: is there a construction analogous to the Bloom-Graham type for germs of CR-regular sets at a singular point?
	A direct generalization turns out to be impossible, since the initial distribution $D_{1}=T_{p}^{c} M_{p}\oplus \overline{T_{p}^{c}} M_{p}$ of complex and anticomplex tangent spaces, which generates the remaining distributions $D_{\nu}$, does not always extend to singular points.
	
	However, we can define the distributions of $D_{\nu}$ at a singular point by considering the sets $D_{\nu}(p)$ of their limit values and obtain a set of invariants composed of the dimensions of the complex linear spans of the sets $D_{\nu}(p)$ (see Section 4). This construction is reminiscent of the construction of the tangent cone to an analytic set and other Whitney cones (see \cite{3}), which serve as analogues of the tangent space at singular points of complex-analytic sets. Our definitions correspond to the definitions of the first, second, and fourth Whitney cones, which seem most natural in this situation.

	Moreover, if we fix one of the three ways of defining the tangent cone, then at the singular point $p$ two sets of invariants arise. These invariants are constructed from the chains $D_{\nu}^{1}(p)$ and $D_{\nu}^{2}(p)$ of linear spaces, which are defined using the linear spans (over $\mathbb{C}$) of the sets $D_{\nu}(p)$ (see Section 4).

	These chains implement different variants of interaction of the real structure of a CR-regular set with the complex structure of the ambient space at a singular point. The first chain is embedded in $\mathbb{C}T_{p} \, M_{p}$ -- an analogue of the complexification of the tangent space, and the second -- in $T_{p}^{c} \, M_{p}\oplus \overline{T_{p}^{c}} \, M_{p}$ -- an analogue of the sum of the complex and anticomplex tangent spaces (for both analogues, the Zariski definition is used, see Section 4).

	This yields a more refined set of invariants compared to the CR-regular point, where exactly one chain is considered (it coincides with the chain of spaces $D_{\nu}^{1}(p)$), and the chain of spaces $D_{\nu}^{2}(p)$ is trivial (all spaces coincide and are equal to $D_{1}^{2}(p)=T_{p}^{c} M_{p}\oplus \overline{T_{p}^{c}} M_{p}$).
	
	We combine the invariants constructed in this way into a single set, which we call the \textit{Bloom-Graham singular type}, or simply the \textit{singular type}.
	
	Note that we potentially have $3 \times 2=6$ invariant sequences, which are obtained by considering three different definitions of the tangent cone and two different tangent spaces at a singular point (real and the sum of complex and anticomplex).

	Moreover, all our constructions can be extended to the case of germs of smooth ($C^{\infty}$) CR-regular sets (i.e., sets locally defined by a system of smooth equations and which are CR manifolds at a generic point). However, to do this, we must additionally impose certain restrictions on the set of singular points: we must require that the set of limit values of the set ${\rm reg} \, M$ coincide with $M$ (so that passage to the limit to singular points along a sequence of non-singular points is possible).
	
	Moreover, one can also consider the weighted Bloom-Graham type (the construction is completely analogous; for the definition of the weighted type, see \cite{4}) and obtain its analogue for the germ of a CR-regular set -- the weighted singular Bloom-Graham type.
	
	We also note that the Bloom-Graham type has two varieties: the geometric type, which we discuss in detail here, and the analytic type, which remains outside the scope of our considerations. The analytic type is defined by a certain standard form of defining relations for a germ of a CR manifold of a fixed geometric type. The Bloom-Graham theorem states the equivalence of both types -- the geometric and the analytic. Perhaps, by analogy with the Bloom-Graham geometric type at a CR-regular point, the singular type contains the key to constructing a certain standard form of equations at an RC-singular point. Moreover, the question arises not only about the analogue of the analytic type, but also about its equivalence to the analogue of the geometric type (i.e., the singular type).
	\vspace{3ex}
	
	The same approach can be applied to constructing the singular $l$-type (see Section 5). The $l$-type is a convenient tool for providing a constructive criterion for the holomorphic non-degeneracy of a CR manifold (see \cite{6}). Without going into details, we note only that holomorphic non-degeneracy is the absence of a non-zero holomorphic vector field tangent to the manifold. Moreover, under the assumption of finite Bloom-Graham type of the germ of a CR manifold, holomorphic non-degeneracy is a criterion for the finite-dimensionality of the algebra of infinitesimal holomorphic automorphisms of this germ (and for a real hypersurface, the finite-type condition need not be additionally imposed -- it is satisfied automatically).
	
	For a non-singular point, the $l$-type is constructed from some expanding chain of vector spaces at a given point $p$: the starting object $\mathcal{L}^{0}(p)$ is the linear span (over $\mathbb{C}$) of the holomorphic gradients of the defining equations, and the $\nu$-space $\mathcal{L}^{\nu}(p)$ is the linear span (over $\mathbb{C}$) of the $\nu$-fold derivations of the said holomorphic gradients by CR fields.
	
	At a singular point, we define the singular $l$-type $\mathfrak{l}(p)$ as the sequence of dimensions of the linear spans (over $\mathbb{C}$) of the set of limit values of the spaces $\mathcal{L}^{\nu}$. The limit values are also defined in three ways and thereby increase the number of invariants, i.e., we potentially have 3 invariant sequences. 
	
	In this way, we use a very fruitful idea, standard in CR geometry, to study singular points -- namely, the interaction between the real and complex structures of a manifold. We obtain sets of invariants at the intersection of CR geometry and singularity theory (see, for example, \cite{7}).

	\vspace{3ex}
	\section{Type in the sense of Bloom and Graham}
	
	\vspace{3ex}
	
	Let $M \subset \mathbb{C}^{r}, \, r=n+k,$ be a real-analytic submanifold defined in a neighborhood $U$ of the origin. Suppose that outside some proper subset $\Sigma \subset M \cap U$ the manifold $M$ is a generic CR manifold of CR-type $(n,k)$ ($n$ is the CR-dimension, $k$ is the codimension). That is, for all points $p\in (M \cap U)\setminus \Sigma$ the dimension of the complex tangent space $T_{p}^{c} M = T_{p} M \cap i T_{p} M$ is constant and equal to $n$, and for points $p\in \Sigma$ either the dimension of the complex tangent space is greater than $n$ (an RC-singular point) or the complex tangent space is undefined (a singular point of a smooth structure).
	
	Note that $\Sigma$ can be the empty set, but this case is not meaningful for us. Clearly, $\Sigma$ is a proper analytic subset of $M \cap U$, and therefore the complement of $\Sigma$ is dense in $M \cap U$.
	
	Points of the set $\Sigma$ are called \textit{singular}. For an RC-singular point, the number ${\rm ord}(p)={\rm dim} \, T_{p}^{c} - n$ is called the \textit{order} of the RC-singular point.
	
	On the set $(M \cap U)\setminus \Sigma$, a biholomorphic invariant -- the Bloom-Graham type -- is defined.
	Recall its definition.
	
	Let $M_{p}$ be the germ of $M$ at the point $p$.
	Denote by $D_{1} = T^{c} \, M\oplus \overline{T^{c}} \, M$ the distribution of complex and anticomplex tangent spaces on $M$ defined in a neighborhood of the point $p \in (M \cap U)\setminus \Sigma$ (i.e., the distribution generated by fields of type $(1,0)$ and $(0,1)$).
	
	Let's define inductively a sequence of distributions
	
	$$D_{\nu+1} = [D_{\nu},D_{1}]+D_{\nu},$$
	where $[\cdot,\cdot]$ is the commutator of vector fields.
	
	Let $D_{\nu}(p)$ be the value of $D_{\nu}$ at $p$. Then we have a sequence of embedded complex linear spaces
	
	$$D_{1}(p) \subseteq D_{2}(p) \subseteq ... \subseteq D_{\nu}(p) \subseteq ... \subseteq \mathcal{D}(p) \subseteq \mathbb{C}T_{p} \, M,$$
	where $\mathcal{D}(p) = \bigcup_{\nu=1}^{\infty} D_{\nu}(p)$, and $\mathbb{C}T_{p} \, M$ is the complexification of the tangent space. Note that the constructed sequence stabilizes after a certain moment.
	
	$M_{p}$ is a germ of finite type if $\mathcal{D}(p)= \mathbb{C}T_{p} \, M$. If $\mathcal{D}(p) \subsetneqq \mathbb{C}T_{p} \, M$, then $M_{p}$ is a germ of infinite type.
	
	Let $d_{\nu}={\rm dim} \, D_{\nu}(p)$.
	
	Let us note all the numbers at which a jump in dimension occurred, i.e., those $\nu$ for which $d_{\nu}>d_{\nu-1}$. We obtain an increasing sequence of natural numbers: $2 \leq m_{1} < ... < m_{l}$. Let $k_{j}$ denote the magnitudes of the jumps, i.e., $k_{j}=d_{m_{j}}-d_{m_{j-1}}$.
	
	The resulting data are combined into a set called the Bloom-Graham type. Namely, if $M_{p}$ is a germ of finite type, then the type is defined as the set
	$m(p) = ((m_{1},k_{1}),...,(m_{l},k_{l}))$; if $M_{p}$ is a germ of infinite type, then the type is defined as the set $m(p) = ((m_{1},k_{1}),...,(m_{l},k_{l}),(\infty,d))$, where $d={\dim} (\mathbb{C}T_{p} \, M)-{\dim} (\mathcal{D}(p))=k-k_{1}-...-k_{l}$. The quantity $d=d(p)$ is called the \textit{defect} of the germ $M_{p}$. For a germ of finite type $d(p)=0$, and for a germ of infinite type $d(p)$ can take values from $1$ to $k$.
	
	\vspace{3ex}
	\section{Whitney cones}
	
	\vspace{3ex}
	
	In paper \cite{5}, six natural definitions of the tangent cone to an analytic set in complex space were introduced (see also \cite{3}). We will consider three of them (the 1st, 2nd, and 4th), which are well suited for our purposes, and give similar, but slightly more general definitions. In doing so, we will consider the following class of sets.

	\vspace{3ex}
	
	\textbf{Definition 1.} By a CR-regular set we mean a locally defined (on some open set $U$) set $M$ of common zeros of real-analytic functions in a complex space such that at a generic point (outside some proper analytic subset $\Sigma$) it is a CR manifold of fixed CR type $(n,k)$. Points at which $M$ is a CR manifold are called CR-regular points. The set of CR-regular points is denoted by ${\rm reg} \, M$.
	
	\vspace{3ex}
	All our considerations and constructions are local and concern germs, and to construct auxiliary objects, we choose some representative of the germ.
	We will consider germs of locally connected (i.e., $M$ must be connected in some neighborhood of a singular point) and irreducible CR-regular sets (considering reducible sets reduces to considering their irreducible components, and these components do not necessarily have the same CR-types and dimensions). Moreover, the singular set $\Sigma\subset M$ can contain both RC-singular points and singular points of a smooth structure.
	\vspace{3ex}
	
	Let $\mathcal{V}$ be the distribution of subspaces of tangent spaces to $M$ defined on ${\rm reg} \, M$, the set of CR-regular points of $M \subset \mathbb{C}^{r}$. Note that the definitions of Whitney cones follow from ours if we set $\mathcal{V}=T M$, replace $M$ with a complex-analytic set, and additionally assume that $V$ is a holomorphic vector field in item 1 of our definition (see below).
	
	\vspace{3ex}
	
	\textbf{Definition 2.} Let $p\in \Sigma$.
	
	(1) (analogue of the 1st cone) A vector $v$ belongs to the set $\mathcal{C}_{1}(M,\mathcal{V},p)$ if there exists an analytic vector field $V$ defined in some neighborhood $U'$ of $p$ and belonging to the distribution $\mathcal{V}$ at each point of ${\rm reg} \, M\cap U'$, such that $V(p)=v$.
	
	(2) (analogue of the 2nd cone) A vector $v$ belongs to the set $\mathcal{C}_{2}(M,\mathcal{V},p)$ if for every $\varepsilon>0$ there exists $\delta>0$ such that if $p'\in {\rm reg} \, M$ and $|p'-p|<\delta$, then $\lVert v'-v \rVert<\varepsilon$ for some $v'\in \mathcal{V}(p')$ (here $|\cdot|$ and $\lVert \cdot \rVert$ denote the standard metrics for finite-dimensional Euclidean spaces for points and vectors, respectively).
	
	(3) (analogue of the 4th cone) Vector $v$ belongs to the set $\mathcal{C}_{3}(M,\mathcal{V},p)$ if there exist sequences of points $p^{j}\in {\rm reg} \, M$ and vectors $v^{j}\in \mathcal{V}(p^{j})$ such that $p^{j}\longrightarrow p$ and $v^{j}\longrightarrow v$ for $j\longrightarrow \infty$.
	
	\vspace{3ex}
	
	It is clear that considering CR-regular sets that are isolated points makes no sense, since the definitions of cones imply the possibility of limiting. Alternatively, instead of the third definition, we can consider a modification in which the sequence of points is replaced by an arbitrary subset $S$ of the CR-regular set under consideration. In our context, it is natural to consider semianalytic sets $S$, i.e., sets defined by a system of analytic equalities and inequalities:
	
	\vspace{3ex}
	
	(3') (analogue of the 4th cone) A vector $v$ belongs to the set $\mathcal{C}_{3}(M,\mathcal{V},p)$ if there exists a semi-analytic set $S$ such that for all $p'\in {\rm reg} \, M\cap S$ such that $p'\longrightarrow p$, one can choose vectors $v'\in \mathcal{V}(p')$ such that $v'\longrightarrow v$.
	
	\vspace{3ex}
	
	\section{Singular type}
	
	\vspace{3ex}
	
	Now we extend the definition given in Section 2 to singular points, i.e., to the set $\Sigma$. To do this, we fix the $j$-th way of defining the tangent cone ($1\leq j \leq 3$). The type consists of two collections (singular types of the 1st and 2nd kinds), which we will construct below. We define each of the collections $\mathfrak{m}_{j}^{q}, \, 1\leq q \leq 2$. That is, potentially at a singular point we obtain $3 \times 2=6$ times more invariants compared to a CR-regular point. To simplify the notation, we will omit the index $j$ (since the construction is the same for all three ways of defining the tangent cone) and write $\mathfrak{m}^{q}$ instead of $\mathfrak{m}_{j}^{q}$.

	The main ingredient of our invariants is the sets of limit values $D_{\nu}(p)=\mathcal{C}_{j}(M,D_{\nu},p)$ of the distributions $D_{\nu}$ on the set $\Sigma$.
	
	The system of invariants will be constructed based on the linear spans of these sets over $\mathbb{C}$ (hereinafter, we will denote the linear span by $\langle \cdot\rangle$).
	
	We will also need the sets:
	
	1) $T_{p} \, M_{p}$, which at a non-smooth point is defined as the tangent Zariski space for the germ $M_{p}$, i.e. as the kernel of the Jacobian matrix. Namely, $T_{p} \, M_{p}$ consists of vectors $a_{1}\frac{\partial}{\partial z_{1}}+...+a_{r}\frac{\partial}{\partial z_{r}}+b_{1}\frac{\partial}{\partial \bar{z}_{1}}+...+b_{r}\frac{\partial}{\partial \bar{z}_{r}},$ such that at point $p$ we have $\sum_{j=1}^{r} (a_{j}\frac{\partial \rho_{s}}{\partial z_{j}}+b_{j}\frac{\partial \rho_{s}}{\partial \bar{z}_{j}})=0$ for all $s, \, 1\leq s \leq k$ (at an RC-singular point the definition is standard). By $\mathbb{C}T_{p} \, M_{p}$ we denote the complexification of the space $T_{p} \, M_{p}$.
	
	2) $T_{p}^{c} \, M_{p}$, which at a nonsmooth point is defined as the complex tangent Zariski space. Namely, $T_{p}^{c} \, M_{p}$ consists of vectors $a_{1}\frac{\partial}{\partial z_{1}}+...+a_{r}\frac{\partial}{\partial z_{r}},$ such that at the point $p$ we have $\sum_{j=1}^{r} a_{j}\frac{\partial \rho_{s}}{\partial z_{j}}=0$ for all $s, \, 1\leq s \leq k$ (at an RC-singular point the definition is standard).
	
	3) $\overline{T_{p}^{c}} \, M_{p}$, which at a nonsmooth point is defined as the anticomplex tangent Zariski space. Namely, $\overline{T_{p}^{c}} \, M_{p}$ consists of vectors $a_{1}\frac{\partial}{\partial \bar{z}_{1}}+...+a_{r}\frac{\partial}{\partial \bar{z}_{r}},$ such that at the point $p$ we have $\sum_{j=1}^{r} a_{j}\frac{\partial \rho_{s}}{\partial \bar{z}_{j}}=0$ for all $s, \, 1\leq s \leq k$ (at an RC-singular point the definition is standard).
	
	4) $\mathcal{D}(p) = \bigcup_{\nu=1}^{\infty} D_{\nu}(p)$.

	\vspace{3ex}
	
	1) \textit{Singular type of the first kind.}
	
	We define $D_{\nu}^{1}(p)$ as the linear span $\langle D_{\nu}(p)\rangle$ of the set $D_{\nu}(p)$.
	We say that $M_{p}$ is a germ of a CR-regular set \textit{of finite singular type of the first kind} if $\mathcal{D}^{1}(p)=\langle\mathcal{D}(p)\rangle=\langle\mathcal{C}_{j}(M,\mathcal{D},p)\rangle$ coincides with $\mathbb{C}T_{p} \, M_{p}$. If $\langle\mathcal{D}(p)\rangle \subsetneqq \mathbb{C}T_{p} \, M_{p}$, then $M_{p}$ is a germ of a CR-regular set \textit{of infinite singular type of the first kind}.
	
	The singular type of the first kind is defined in exactly the same way as the Bloom-Graham type at a non-singular point: here the spaces $D_{\nu}^{1}(p)=\langle D_{\nu}(p)\rangle$ act as a chain of distributions. But there is one exception. If the dimension ${\rm dim}_{\mathbb{C}} \, D_{1}^{1}(p)$ exceeds twice the CR-dimension $2 n$ at a generic point (this means that $D_{1}$ does not extend to $p$), then we set $m_{1}=1, \, k_{1}={\rm dim}_{\mathbb{C}} \, D_{1}^{1}(p)-2 n$. In this case, the remaining values $(m_{j},k_{j})$ for $j\geq 2$ fix the number of the distribution on which the jump in dimension occurs and the magnitude of this jump.
	
	If ${\rm dim}_{\mathbb{C}} \, D_{1}^{1}(p) = 2 n$, then $m_{1}\geq 2$, and all quantities $(m_{j},k_{j})$ are defined literally in the same way as at a non-singular point. Note that in this case, if $D_{1}$ extends to $p$, then the spaces $D_{\nu}$ can be defined in the same way as above using commutators, without using the limit passage (both definitions will be equivalent).
	
	Thus, at a singular point we can define a set $\mathfrak{m}^{1}(p) = ((m_{1},k_{1}),...,(m_{l},k_{l}))$ (if the type is finite) or a set $\mathfrak{m}^{1}(p) = ((m_{1},k_{1}),...,(m_{l},k_{l}),(\infty,d^{1}))$ (if the type is infinite), where $d^{1}=d^{1}(p)={\dim} (\mathbb{C}T_{p} \, M)-{\dim} (\mathcal{D}^{1}(p))$. We call the set $\mathfrak{m}^{1}$ the \textit{singular type of the first kind}, and the quantity $d^{1}(p)$ the \textit{singular defect of the first kind}.
	
	\vspace{3ex}

	Note that at a non-singular point the type of the first kind coincides with the usual Bloom-Graham type.
	
	\vspace{3ex}
	
	2) \textit{Singular type of the 2nd kind.}
	
	Next, at a singular point, we can define more subtle invariants.
	Namely,
	consider the sequence of nested linear spaces $D_{\nu}^{2}(p)= \langle D_{\nu}(p)\rangle\cap (T_{p}^{c} \, M_{p}\oplus \overline{T_{p}^{c}} \, M_{p})$:

	\begin{equation}\label{eq1}
	D_{1}^{2}(p) \subseteq D_{2}^{2}(p) \subseteq ... \subseteq D_{\nu}^{2}(p) \subseteq ... \subseteq \mathcal{D}^{2}(p) \subseteq  T_{p}^{c} \, M_{p}\oplus \overline{T_{p}^{c}} \, M_{p},
	\end{equation}	
	where $\mathcal{D}^{2}(p) = \bigcup_{\nu=1}^{\infty} D_{\nu}^{2}(p)$. Note that the constructed sequence stabilizes after a certain moment. That is, we consider a sequence of nested subspaces inside the space $T_{p}^{c} \, M_{p}\oplus \overline{T_{p}^{c}} \, M_{p}$.
	
	We now repeat the construction described in Section 2 for the sequence of subspaces \eqref{eq1}. For simplicity, we use the same notation for weights and multiplicities as in the definition of the type of the first kind -- $m_{j}$ and $k_{j}$ (however, these quantities may be different).

	First, we say that $M_{p}$ is a germ of a CR-regular set \textit{of finite singular type of the 2nd kind} if $\mathcal{D}^{2}(p)$ coincides with $T_{p}^{c} \, M_{p}\oplus \overline{T_{p}^{c}} \, M_{p}$. If $\mathcal{D}^{2}(p) \subsetneqq T_{p}^{c} \, M_{p}\oplus \overline{T_{p}^{c}} \, M_{p}$, then $M_{p}$ is a germ of a CR-regular set \textit{of infinite singular type of the 2nd kind}.
	
	Let $d_{\nu}={\rm dim}_{\mathbb{C}} \, D_{\nu}^{2}(p)$.
	
	Let's mark all the numbers where the dimension jumps, i.e., those $\nu$ for which $d_{\nu}>d_{\nu-1}$. We obtain an increasing sequence of natural numbers: $1 \leq m_{1} < ... < m_{l}$.
	Moreover, if $m_{1}=1$, this means that the distribution $D_{1}$ does not extend to the point $p$. Let $k_{j}$ denote the magnitudes of the jumps, i.e., $k_{j}=d_{m_{j}}-d_{m_{j-1}}$.

	The obtained data are combined into a set, which we will call the \textit{singular type of the second kind}. Namely, if $M_{p}$ is a germ of a CR-regular set of finite singular type of the second kind, then the singular type of the second kind is defined as the set
	$\mathfrak{m}^{2}(p) = ((m_{1},k_{1}),...,(m_{l},k_{l}))$; if $M_{p}$ is a germ of a CR-regular set of infinite singular type of the second kind at $p$, then the singular type of the second kind is defined as the set $\mathfrak{m}^{2}(p) = ((m_{1},k_{1}),...,(m_{l},k_{l}),(\infty,d^{2}))$, where $d^{2}={\dim} (T_{p}^{c} \, M_{p}\oplus \overline{T_{p}^{c}} \, M_{p})-{\dim} \langle\mathcal{D}^{2}(p)\rangle$. The quantity $d^{2}=d^{2}(p)$ will be called the \textit{singular defect of the 2nd kind} of the germ $M_{p}$ of a CR-regular set $M$. For a germ of a CR-regular set of finite singular type of the 2nd kind, $d^{2}(p)=0$, and for a germ of a CR-regular set of infinite type, $d^{2}(p)$ can take values from $1$ to ${\dim} (T_{p}^{c} \, M_{p}\oplus \overline{T_{p}^{c}} \, M_{p})$.
	
	\vspace{3ex}
	
	Thus, we obtain the sets $\mathfrak{m}^{q}(p)$ -- \textit{singular type of the $q$-th kind} and the value $d^{q}(p)$ -- \textit{singular defect of the $q$-th kind} for $1\leq q \leq 2$.
	
	To construct similar sets for germs of reducible sets, we need to replace the sets $D_{\nu}^{1}(p)$ and $D_{\nu}^{2}(p)$ with the unions of all such sets for all irreducible components.

	\vspace{3ex}

	Thus, we obtain an invariant set $(\mathfrak{m}^{1},\mathfrak{m}^{2})$, which we will call the \textit{singular type}, and an invariant set $(d^{1},d^{2})$, which we will call the \textit{singular defect}. The biholomorphic invariance of both sets is clear from the definition, and we obtain the following statement.
	
	\vspace{3ex}
	
	\textbf{Theorem 3.} The singular type and singular defect are biholomorphic invariants of the germ of a CR-regular set.

	\vspace{3ex}

	\textbf{Remark 4.} As a starting invariant at a CR-regular point, we can also consider the weighted Bloom-Graham type (the construction is completely analogous; for the definition of the weighted type, see \cite{4}). The difference is that instead of distributions $D_{\nu}$, we should consider their weighted analogues $D_{\nu}(\mu)$ for a chosen system of weights $\mu=(\mu_{1},...,\mu_{n})$ of coordinates $z=(z_{1},...,z_{n})$ of the complex tangent space. We thus obtain an analogue of the weighted Bloom-Graham type for a CR-regular set -- the weighted Bloom-Graham singular type $\mathfrak{m}(\mu)$.
	
	\vspace{3ex}
	\section{Singular $l$-type}
	
	\vspace{3ex}
	
	Let us recall the definition of $l$-type.
	
	\vspace{3ex}
	
	\textbf{Definition 5.} Let $z=(z_{1},...,z_{r})$ be coordinates in $\mathbb{C}^{r}$. A \textit{CR field} on a CR manifold $M$ is a vector field tangent to $M$ of the form
	
	$$F_{1}(z,\bar{z})\frac{\partial}{\partial \bar{z}_{1}}+...+F_{r}(z,\bar{z})\frac{\partial}{\partial \bar{z}_{r}},$$
	
	where $F_{j}$ are smooth complex-valued functions. Note that here we are considering fields of type $(0,1)$.
	
	\vspace{3ex}
	
	Let $M$ be a CR manifold of CR type $(n,k)$ defined in a neighborhood of a point $p$ by a system of defining equations $\{\rho_{j}=0, \, 1\leq j \leq k\}$, $L_{1},...,L_{n}$ be a basis of CR vector fields in a neighborhood of a point $p$ on $M$, $\frac{\partial \rho_{j}}{\partial z}=(\frac{\partial \rho_{j}}{\partial z_{1}},...,\frac{\partial \rho_{j}}{\partial z_{r}})$, $L^{\alpha}=L_{1}^{\alpha_{1}}\cdot...\cdot L_{r}^{\alpha_{r}}$. 
	
	\vspace{3ex}
	
	\textbf{Definition 6.} A germ $M_{p}$ of a manifold $M$ is called  \textit{$l$-nondegenerate} if the complex linear span of the system of vectors $\{L^{\alpha}\Big(\frac{\partial \rho_{j}}{\partial z}\Big)(p,\bar{p}), \, |\alpha|\leq l, \, 1\leq j \leq k\}$ coincides with $\mathbb{C}^{r}$.
	
	\vspace{3ex}
	
	\textbf{Definition 7.} The $l$-type of the germ $M_{p}$ of the CR manifold $M$ is the set of dimensions of the complex linear spans of vectors $\mathcal{L}^{\nu}(p)=\{L^{\alpha}\Big(\frac{\partial \rho_{j}}{\partial z}\Big)(p,\bar{p}), \, |\alpha|\leq \nu, \, 1\leq j \leq k\}, \, 1\leq \nu \leq l$, i.e., the set $l(p)=({\rm dim} \, \mathcal{L}^{1}(p),{\rm dim} \, \mathcal{L}^{2}(p),...,{\rm dim} \, \mathcal{L}^{l}(p))$.
	
	\vspace{3ex}
	
	At a singular point, we define the singular $l$-type $\mathfrak{l}(p)$ as the sequence of dimensions of the linear spans (over $\mathbb{C}$) of the set of limit values of the spaces $\mathcal{L}^{\nu}(p), \nu\geq 0$. The limit values can also be defined in three ways by analogy with the cones $\mathcal{C}_{j}, \, 1\leq j \leq 3$. The difference from the above definitions of the cones $\mathcal{C}_{j}(M,\mathcal{V},p)$ is that instead of the distribution $\mathcal{V}$, we need to consider a set of linear spaces $\mathcal{L}^{\nu}(p)$ that depend analytically on the point. The spaces $\mathcal{L}^{\nu}(p)$ are not subspaces of the complexification of the tangent space $\mathbb{C} T_{p} \, M_{p}$, but for our purposes this is unimportant. Thus, we obtain the following set of values: $\mathfrak{l}(p)=({\rm dim} \, \langle\mathcal{C}_{j}(M,\mathcal{L}^{0},p)\rangle ,{\rm dim} \, \langle\mathcal{C}_{j}(M,\mathcal{L}^{1},p)\rangle,...,{\rm dim} \, \langle\mathcal{C}_{j}(M,\mathcal{L}^{l},p)\rangle)$ for each $j, \, 1\leq j\leq 3$. Since the limit values are defined in three ways, at a singular point we potentially have three times more invariants than at a CR-regular point. Unlike the type at a regular point, $\nu$ can take the value zero (at a regular point we always have ${\rm dim} \, \mathcal{L}^{0}(p)=k$). The definition of $l$-nondegeneracy also carries over to the case of a singular point: a germ $M_{p}$ of a CR-regular set $M$ is called \textit{$l$-nondegenerate} if the complex linear span of the set of limit values of the space $\mathcal{L}^{l}(p)$ coincides with $\mathbb{C}^{r}$.
	We have the following statement.
	
	\vspace{3ex}
	
	\textbf{Theorem 8.} The singular $l$-type is a biholomorphic invariant of the germ of a CR-regular set.
	
	\vspace{3ex}
	
	\textbf{Remark 9.} At a CR-regular point, the singular $l$-type coincides with the $l$-type.

	\section{Examples}
	\vspace{3ex}
	Let us illustrate these definitions with a series of examples. Their main purpose is to present examples of germs of CR-regular sets for which the type of the second kind is nontrivial. It is clear that all possible examples of types of the first kind, as well as examples of singular $l$-types, are realized using germs of CR manifolds.

	\vspace{3ex}
	In the following example, the distribution $D_{1}^{q}$ does not extend to an RC-singular point:
	
	\vspace{3ex}
	
	\textbf{Example 10.} In this example, we will use the third definition of the tangent cone. Let $(z_{1},z_{2},z_{3})$ be coordinates in $\mathbb{C}^{3}$. Consider a CR-regular set $M$ defined by the system
	
	$$M=\{|z_{1}|^{2}+|z_{2}|^{2}+|z_{3}|^{2}=1, \, {\rm Im} \, z_{1}=0\}.$$
	This is a section of a 5-dimensional sphere by a hyperplane, i.e., a 4-dimensional sphere.
	
	To calculate singular types, consider the field $X=a_{1}\frac{\partial}{\partial z_{1}}+a_{2}\frac{\partial}{\partial z_{2}}+a_{3}\frac{\partial}{\partial z_{3}}$.
	Let us write down the conditions of tangency of the field $X$ of the CR-regular set $M$:
	
	$$a_{1}\bar{z}_{1}+a_{2}\bar{z}_{2}+a_{3}\bar{z}_{3}=0, \, \, \, \frac{a_{1}}{2 i}=0.$$
	We have $a_{1}=0$, at the points $(\pm 1,0,0)$ the dimension of the complex tangent space is two (since $a_{2}$ and $a_{3}$ can take arbitrary values), and at the other points $a_{2}\bar{z}_{2}=-a_{3}\bar{z}_{3}$, and the dimension of the complex tangent space is one (since at least one of the quantities $\bar{z}_{2}$ and $\bar{z}_{3}$ does not vanish). Therefore, $(\pm 1,0,0)$ are RC-singular points (the singular set is zero-dimensional). It is clear that along the direction $z_{3}=\alpha z_{2}$ the quantity $a_{2}$ will have the limit $-a_{3}\bar{\alpha}$ for $z_{2},z_{3}\longrightarrow 0$. Therefore the field $X$ does not extend to RC-singular points. Moreover, all possible pairs $(a_{2},a_{3})$ of limit values of the coefficients can be arbitrary. And also at RC-singular points $T_{p}^{c} \, M_{p}\oplus \overline{T_{p}^{c}} \, M_{p}=\mathbb{C}T_{p} \, M_{p}$. From this we obtain that types of the 1st and 2nd kind coincide and are equal to $((1,2))$ at all RC-singular points.
	
	\vspace{3ex}
	And in the next series of examples, the distribution $D_{1}^{q}$ extends to RC-singular points:
	
	\vspace{3ex}
	\textbf{Example 11.} 1) In this example, we will use the first definition of the tangent cone. Let $(z_{1},z_{2},z_{3})$ be coordinates in $\mathbb{C}^{3}$. Consider a CR-regular set $M$ defined by the system
	
	$$M=\{|z_{1}|^{2}+|z_{2}|^{2}(1+|z_{3}|^{2})=1, \, {\rm Im} \, z_{1}=0\}.$$
	At a generic point, this is a 4-dimensional irreducible non-compact (since for $z_{2}=0$ the variable $z_{3}$ can take any values) manifold.
	
	Let $X=a_{1}\frac{\partial}{\partial z_{1}}+a_{2}\frac{\partial}{\partial z_{2}}+a_{3}\frac{\partial}{\partial z_{3}}$. Let us write down the tangency conditions of the field $X$ of the CR-regular set $M$:
	
	$$a_{1}\bar{z}_{1}+a_{2}\bar{z}_{2}(1+|z_{3}|^{2})+a_{3}|z_{2}|^{2}\bar{z}_{3}=0, \, \, \, a_{1}=0.$$
	We have $a_{1}=0$, the points $(\pm 1,0,z_{3})$ are RC-singular for all $z_{3}$ (the dimension of the complex tangent space is two; the singular set is real two-dimensional), and at the remaining points $a_{2}=-a_{3}\frac{z_{2}\bar{z}_{3}}{1+|z_{3}|^{2}}$ (the dimension of the complex tangent space is one). Moreover, the coefficient at $\frac{\partial}{\partial z_{2}}$ of the field $X$ is divisible by $z_{2}$, which means that the dimensions of the spaces $D_{\nu}^{q}$ at RC-singular points cannot increase. 
	And also at RC-singular points $T_{p}^{c} \, M_{p}\oplus \overline{T_{p}^{c}} \, M_{p}=\mathbb{C}T_{p} \, M_{p}$, and in some neighborhood of each RC-singular point $M$ is connected. From this we obtain that the types of the 1st and 2nd kind coincide and are equal to $((\infty,2))$ at all RC-singular points. 
	
	\vspace{3ex}

	2) In this example, we will use the first definition of the tangent cone. Let $(z_{1},z_{2},z_{3},z_{4})$ be coordinates in $\mathbb{C}^{4}$. Consider a CR-regular set $M$ defined by the system
	
	$$M=\{|z_{1}|^{2}+|z_{2}|^{2}=1, \, {\rm Im} \, z_{1}=0, \, {\rm Im} \, z_{3}=|z_{4}|^{2}\}.$$
	This is the direct product of a 2-dimensional sphere and a 3-dimensional hyperquadric, which is projectively equivalent to the 3-dimensional sphere. That is, $M$ is a 5-dimensional connected irreducible non-compact (since $z_{3}$ and $z_{4}$ can take arbitrarily large absolute values) manifold.
	
	Let $X=a_{1}\frac{\partial}{\partial z_{1}}+a_{2}\frac{\partial}{\partial z_{2}}+a_{3}\frac{\partial}{\partial z_{3}}+a_{4}\frac{\partial}{\partial z_{4}}$. Let us write down the tangency conditions of the field $X$ of a CR-regular set $M$:
	
	$$a_{1}\bar{z}_{1}+a_{2}\bar{z}_{2}=0, \, \, \, a_{1}=0, \, \, \, \frac{a_{3}}{2 i}=a_{4}\bar{z}_{4}.$$
	We have $a_{1}=0, \, a_{3}=2 i a_{4}\bar{z}_{4}$, the points $(\pm 1,0,z_{3},z_{4})$ are RC-singular for all $z_{3}, z_{4}$ (the dimension of the complex tangent space is two), and at the remaining points $a_{2}=0$ and $a_{4}$ is a free parameter (the dimension of the complex tangent space is one). The field $X$ extends to RC-singular points, and the commutator of depth one gives a new direction (a non-zero coefficient at $\frac{\partial}{\partial z_{3}}$ and $\frac{\partial}{\partial \bar{z}_{3}}$). But the directions $\frac{\partial}{\partial z_{2}}$ and $\frac{\partial}{\partial \bar{z}_{2}}$ cannot be obtained, since at a generic point the coefficient $a_{2}$ is equal to zero. We have: $\mathfrak{m}^{1}=((2,1),(\infty,2)), \, \mathfrak{m}^{2}=((\infty,2))$ at all RC-singular points.

	\vspace{3ex}

	3) The following example demonstrates that the weight $m_{1}$ in a type of the 2nd kind can take arbitrary values greater than 2. In this example, we will use the first definition of the tangent cone. Let $(z_{1},...,z_{7})$ be coordinates in $\mathbb{C}^{7}$. Consider a CR-regular set $M$ defined by the system
	
	$${\rm Im} \, z_{1}=0, \,  |z_{1}|^{2}+|z_{2}|^{2}(1+(z_{3}^s+\bar{z}_{3}^s)|z_{3}|^{2})=1,$$
	
	$$|z_{1}|^{2}+|z_{6}|^{2}+|z_{4}|^{2}(1+(z_3^s+\bar{z}_{3}^s)|z_{3}|^{2})=1,$$
	
	$$|z_{1}|^{2}+|z_{7}|^{2}+(z_5+\bar{z}_{5})|z_{2}|^{2}|z_{4}|^{2}(1+(z_3^s+\bar{z}_{3}^s)|z_{3}|^{2})=1,$$
	where $s$ is a non-negative integer.
	At a generic point, this is a 10-dimensional irreducible non-compact manifold (since for $z_{2}=z_{4}=0$, the variable $z_{3}$ can take any value).
	
	Let $X=a_{1}\frac{\partial}{\partial z_{1}}+...+a_{7}\frac{\partial}{\partial z_{7}}$. We write down the tangency conditions of the field $X$ of a CR-regular set $M$:
	
	$$\frac{a_{1}}{2 i}=0, \,	\, \,
	a_1 \bar{z}_{1}+a_2 \bar{z}_{2}(1+(z_3^s+\bar{z}_{3}^s)|z_{3}|^{2})+a_3 |z_{2}|^{2}(s z_3^s \bar{z}_{3}+(z_3^s+\bar{z}_{3}^s)\bar{z}_{3})=0,$$
	
	$$a_1 \bar{z}_{1}+a_3 |z_{4}|^{2}(s z_3^s\bar{z}_{3}+(z_3^s+\bar{z}_{3}^s)\bar{z}_{3})+a_4 \bar{z}_{4}(1+(z_3^s+\bar{z}_{3}^s)|z_{3}|^{2})+a_6 \bar{z}_{6}=0,$$
	
	$$a_1 \bar{z}_{1}+a_2(z_5+\bar{z}_{5})\bar{z}_{2}|z_{4}|^{2}(1+(z_3^s+\bar{z}_{3}^s)|z_{3}|^{2})+$$
	
	$$+a_3(z_5+\bar{z}_{5})|z_{2}|^{2}|z_{4}|^{2}(s z_3^s\bar{z}_{3}+(z_3^s+\bar{z}_{3}^s)\bar{z}_{3})+a_4(z_5+\bar{z}_{5})|z_{2}|^{2}\bar{z}_{4}(1+(z_3^s+\bar{z}_{3}^s)|z_{3}|^{2})+$$
	
	$$+a_5|z_{2}|^{2}|z_{4}|^{2}(1+(z_3^s+\bar{z}_{3}^s)|z_{3}|^{2})+a_7 \bar{z}_{7}=0.$$
	
	We have $a_{1}=0$, and at a generic point defined by the conditions $\{(|z_{3}|^{2}z_3^s+|z_{3}|^{2}\bar{z}_{3}^s+1)\neq 0, \, z_{2}\neq 0, \, z_{4}\neq 0\}$, the coefficients $a_{2},a_{4},a_{5}$ are expressed through $a_{3},a_{6},a_{7}$ as follows:
	
	$$a_2 = \frac{-a_3\bar{z}_{3}z_2(s z_3^s+z_3^s+\bar{z}_{3}^s)}{(|z_{3}|^{2}z_3^s+|z_{3}|^{2}\bar{z}_{3}^s+1)},$$
	
	$$a_4 = \frac{-a_3\bar{z}_{3}z_4(s z_3^s+z_3^s+\bar{z}_{3}^s)\bar{z}_{4}-a_{6}\bar{z}_{6}}{\bar{z}_{4}(|z_{3}|^{2}z_3^s+|z_{3}|^{2}\bar{z}_{3}^s+1)},$$
	
	$$a_5 = \frac{|z_{2}z_{4}|^{2}\bar{z}_{3}a_3(s z_3^s+z_3^s +\bar{z}_{3}^s)(\bar{z}_{5}+ z_5)+a_6 \bar{z}_{2}\bar{z}_{5}\bar{z}_{6}z_2+a_6\bar{z}_{2}\bar{z}_{6}z_2 z_5+a_{7}\bar{z}_{7}}{|z_{2}z_{4}|^{2}(|z_{3}|^{2}z_3^s+|z_{3}|^{2}\bar{z}_{3}^s+1)}.$$
	Thus, at a generic point, the dimension of the complex tangent space is three.
	The points $(\pm 1,0,z_{3},0,z_{5},0,0)$ are singular points of the smooth structure (the dimension of the space $T_{p}^{c} \, M_{p}\oplus \overline{T_{p}^{c}} \, M_{p}$ is twelve), since for $z_{2}=z_{4}=z_{6}=z_{7}=0$ the coefficients $a_{2},a_{3},a_{4},a_{5},a_{6},a_{7}$ can take arbitrary values.
	
	Consider the singular point $p=(1,0,0,0,1,0,0)$. From the formula for $a_{4}$, it is clear that for $a_{6}\neq 0$, the field $X$ does not extend to the point $p$, therefore $a_{6}=0$. And from the formula for $a_{5}$, it is clear that for $a_{6}=0,a_{7}\neq 0$, the field also does not extend to the point $p$, therefore $a_{7}=0$. For $a_{6}=a_{7}=0$, the formulas for $a_{4}$ and $a_{5}$ take a simpler form:
	
	$$a_4 = \frac{-a_3\bar{z}_{3}z_4(s z_3^s+z_3^s+\bar{z}_{3}^s)}{(|z_{3}|^{2}z_3^s+|z_{3}|^{2}\bar{z}_{3}^s+1)},$$
	
	$$a_5 = \frac{\bar{z}_{3}a_3(s z_3^s+z_3^s +\bar{z}_{3}^s)(\bar{z}_{5}+ z_5)}{(|z_{3}|^{2}z_3^s+|z_{3}|^{2}\bar{z}_{3}^s+1)}.$$
	
	
	
	We define the fields $X_{\nu}$ inductively as follows: $X_{1}=[\bar{X},X], \, X_{\nu}=[\bar{X},X_{\nu-1}]$. The field $X_{\nu}$ belongs to the distribution $D_{\nu+1}$. Moreover, the field $X$ contains the term $a_3\bar{z}_{3}^{s+1}(\bar{z}_{5}+ z_5)\frac{\partial}{\partial z_{5}}$ (if we expand the expression for $a_{5}$ in a Taylor series), and the field $\bar{X}$ contains the term $\bar{a}_{3}\frac{\partial}{\partial \bar{z}_{3}}$. Therefore, the field $X_{s+1}$ contains the term $c \frac{\partial}{\partial z_{5}}$ with a non-zero coefficient $c=(s+1)! a_3\bar{a}_{3}^{s+1}(\bar{z}_{5}+ z_5)$. Moreover, at the point $p$ it is impossible to obtain a non-zero coefficient at $\frac{\partial}{\partial z_{5}}$ in a number of commutations less than $s+1$, since the degree in $z_{3},\bar{z}_{3}$ of the coefficient at $\frac{\partial}{\partial z_{5}}$ of the field $X$ is equal to $s+1$. And the coefficients at $\frac{\partial}{\partial z_{2}}, \frac{\partial}{\partial \bar{z}_{2}}, \frac{\partial}{\partial z_{4}}, \frac{\partial}{\partial \bar{z}_{4}}$ of the fields $X,\bar{X}$ are multiplied by $z_{2}, \bar{z}_{2}, z_{4}, \bar{z}_{4}$ respectively, therefore fields from $D_{\nu}^{2}(p)$ can contain only terms $\frac{\partial}{\partial z_{3}}, \frac{\partial}{\partial \bar{z}_{3}}, \frac{\partial}{\partial z_{5}}, \frac{\partial}{\partial \bar{z}_{5}}$ with some coefficients. Also at the point $p$ we have $T_{p}^{c} \, M_{p}\oplus \overline{T_{p}^{c}} \, M_{p}=\mathbb{C}T_{p} \, M_{p}$, and in some neighborhood of each singular point $M$ is connected. Moreover, at the generic point, the coordinates in the complex tangent space are $z_{3},z_{6},z_{7}$, therefore non-zero terms of the form $c_{1}\frac{\partial}{\partial z_{5}}, c_{2}\frac{\partial}{\partial \bar{z}_{5}}$ for arbitrary $c_{1}, c_{2}\in\mathbb{C},$ cannot be obtained in fields from $\mathcal{D}^{2}(p)$ (otherwise $z_{5}$ would also be a coordinate in the complex tangent space). This means that the coefficients $c_{1},c_{2}$ are dependent, i.e. the magnitude of the jump in the dimension of the space $D_{s+1}^{2}(p)$ is equal to unity.
	
	Therefore we have:  $\mathfrak{m}^{1}=\mathfrak{m}^{2}=((s+2,1),(\infty,9))$ at the point $p$.

	\vspace{3ex}
	4) Example 3) allows us to demonstrate that the weights $m_{j}$ in a type of the second kind can take arbitrary values greater than two, and with arbitrary multiplicities $k_{j}$. To do this, we need to consider the direct products of the germs from Example 3) for various $s$. By choosing suitable values of $s$, we can obtain arbitrary types of the first and second kind up to some defect $d$.

	\vspace{3ex}
	\textbf{Example 12.} In this example, we will use the third definition of the tangent cone. Let $(z_{1}=x_{1}+i y_{1},z_{2}=x_{2}+i y_{2},z_{3}=x_{3}+i y_{3})$ be coordinates in $\mathbb{C}^{3}$. Consider the CR-regular set $M$ (quadratic cone) given by the equation
	
	$$x_{1}^{2}=x_{2}^{2}+x_{3}^{2}.$$
	
	The origin $p$ is a nonsmooth point, at which $T_{p}^{c} \, M_{p}\oplus \overline{T_{p}^{c}} \, M_{p}=T_{p} \, M_{p}=\mathbb{C}^{3}$.
	
	The singular types of the first and second kind coincide and are equal to $((1,2))$, since at the generic point the dimension of the complex tangent space is two, and the distribution $T_{p}^{c} \, M_{p}\oplus \overline{T_{p}^{c}} \, M_{p}$ does not extend to $p$. The singular $l$-type is $(3)$, since the linear span of the set of limit values of the distribution generated by the holomorphic gradient coincides with $\mathbb{C}^{3}$.

	\vspace{3ex}
	The Bloom-Graham theorem allows us not only to write out a set of invariants but also to specify a standard form to which the defining equations of a germ with a given Bloom-Graham type can be reduced. The reverse approach is also possible: given a standard form, the Bloom-Graham type can be uniquely reconstructed. The same question arises for the RC type. There are some results in this direction (constructing a standard form at an RC singular point) (see, for example, \cite{2}), but a complete picture is lacking. Moreover, in \cite{2}, a normal form was obtained for the case $n=0$ (at a generic point),
	which we did not consider.
	
	\vspace{3ex}
	
	\textbf{Question 13.} To what form can the defining equations of the germ of a CR-regular set of a given singular type at an RC-singular point be reduced? (It is natural to call this form \textit{singular analytic type} by analogy with the case of a CR-regular point.) Does an analogue of the Bloom-Graham theorem on the equivalence of analytic and geometric types hold? (That is, are the singular analytic type and the singular type that we constructed above equivalent?)

	\vspace{3ex}
	
	It is clear, however, that the sets of numbers in types of different kinds are not independent. There are obvious restrictions (at the level of linear algebra) arising from the fact that types of different kinds are constructed from the same collection of sets $D_{\nu}(p)$.
	
	For example, the sum of the multiplicities $k_{1}^{q}$ of the weights $m_{1}^{q}$ in a type of kind $q$ such that $m_{1}^{q}\leq \mu$ for some $\mu$ is not greater than the sum of the multiplicities $k_{2}^{q}$ of the weights $m_{2}^{q}$ in a type of kind $q$ such that $m_{2}^{j}\leq \mu$ (since the cone $\mathcal{C}_{1}$ is embedded in the cone $\mathcal{C}_{2}$, as is clear from the definition). The same is true for pairs of multiplicities $k_{2}^{q},k_{3}^{q}$ (since the cone $\mathcal{C}_{2}$ is embedded in the cone $\mathcal{C}_{3}$), for pairs of multiplicities $k_{j}^{2},k_{j}^{1}$ (since the space $D_{\nu}^{2}(p)$ is embedded in the space $D_{\nu}^{1}(p)$).
	
	\vspace{3ex}
	
	\textbf{Question 14.} Are there any other restrictions on singular types of different kinds? That is, is it possible for any singular type (with the specified restrictions) to produce a germ of a CR-regular set at a singular point with that type? How are weights and multiplicities related in types of different kinds (including for different methods of transition to limit values)?
	
	\vspace{3ex}
	
	The question of the structure of sets of singular type and singular $l$-type also arises. Namely, for CR-regular points, sets on which both types are constant are semianalytic. Therefore, for singular types, the same answer can be expected under certain restrictions.
	
	\vspace{3ex}
	
	\textbf{Question 15.} Under what restrictions on a CR-regular set are those sets on which the Bloom-Graham singular type or singular $l$-type is constant semi-analytic?

	\renewcommand{\refname}{References}

\end{document}